\documentclass[leqno,12pt]{article} 
\usepackage{amsmath, amssymb}
\usepackage{amsthm} 
\theoremstyle{plain} 
\newtheorem{theorem}{\indent\sc Theorem}[section]
\newtheorem{lemma}[theorem]{\indent\sc Lemma}
\newtheorem{corollary}[theorem]{\indent\sc Corollary}
\newtheorem{proposition}[theorem]{\indent\sc Proposition}

\newtheorem{problem}[theorem]{\indent\sc Problem}
\theoremstyle{definition} 
\newtheorem{definition}[theorem]{\indent\sc Definition}

\newtheorem{example}[theorem]{\indent\sc Example}

\newcommand{\zz}{\mathbb{Z}}

\newcommand{\rr}{\mathbb{R}}

\newcommand{\xx}{\mathbf{x}}
\newcommand{\yy}{\mathbf{y}}
\newcommand{\vv}{\mathbf{v}}
\newcommand{\uu}{\mathbf{u}}

\newcommand{\conv}{\mathrm{conv}}
\newcommand{\cut}{{\mathrm{Cut}^{\square}}}

\makeatletter
\def\address#1#2{\begingroup
\noindent\parbox[t]{7.8cm}{%
\small{\scshape\ignorespaces#1}\par\vskip1ex
\noindent\small{\itshape E-mail address}%
\/: #2\par\vskip4ex}\hfill%
\endgroup}%
\makeatother
\title{\uppercase{Compressed Polytopes and Statistical Disclosure Limitation}} 
\author{
\textsc{Seth Sullivant$^*$} 
}
\date{} 
\begin{document}

\maketitle

\footnote{ 
$^{*}$Supported by a NSF graduate research fellowship.
}

\begin{abstract}
We provide a characterization of the compressed lattice polytopes in
terms of their facet defining inequalities and  we
show that every compressed lattice polytope is affinely isomorphic
to a $0/1$-polytope.   As an application, we characterize those graphs
whose cut  polytopes are compressed and discuss
consequences for studying linear programming relaxations in
statistical disclosure limitation.
\end{abstract}

\section{Introduction}

A lattice polytope $P$ is called compressed if every pulling
triangulation of $P$ using only the lattice points in $P$ is unimodular.
Compressed polytopes are  natural  
to study because they represent a more inclusive class of
polytopes than the unimodular polytopes (polytopes where every
triangulation is unimodular).  Furthermore, many naturally occurring  
polytopes are compressed.  An important example is the
Birkhoff polytope of doubly stochastic matrices as shown in \cite{St}.
In fact, the compressed nature of the Birkhoff polytope played a
crucial role in the work of Diaconis and Sturmfels \cite{DS} for the
statistical analysis of ranked data.  
Ohsugi and Hibi's paper \cite{OH} contains many other 
examples.  In this paper, we characterize the compressed polytopes by
their facet defining inequalities, extending a result from \cite{OH}. 

Part of our motivation for studying compressed polytopes comes from
their appearance in algebraic statistics: the marginal polytopes of
decomposable hierarchical models are compressed.  Due to the presence  of a
transitive symmetry group on these marginal polytopes, the connections
between compressed polytopes and 
certain optimization problems in statistical disclosure limitation
are quite deep.  As an application of our main result on compressed
polytopes, we will show that the linear programming relaxations for
maximizing cell entries given marginal sums yield sharp integer bounds
for all values of the 
marginals if and only if the marginal polytope $P_\Delta$ is compressed.  
Coupled with some results about compressed cut polytopes, we are able
to describe some new nondecomposable families of marginals where the linear
programming relaxation yields sharp integer bounds for the
maximization problems.

Here is the outline for our paper.  In the next section we prove the
main result classifying compressed polytopes by their facet defining
inequalities.  We also show that every compressed polytope is affinely
isomorphic to a 0/1 polytope and prove a result about pulling
triangulations for highly symmetric polytopes.  In Section 3, we apply
the main result to characterize the compressed cut polytopes.  
In Section 4 we explain the connection between compressed polytopes
and linear optimization.  Section
5 is devoted to applications of our results in statistical disclosure
limitation which provides new families of marginals where linear
programming yields sharp upper bounds on cell entries.  These results
also suggest families in which to search for large integer
programming gaps \cite{HS}.  

\section{Characterization of compressed polytopes}

In this section, we derive our main result about the structure of the
facet definining inequalities of compressed polytopes.  We assume the
reader is familiar with polyhedral geometry and regular subdivisions.
A standard reference for this material is \cite{Z}.

\begin{definition}
Let $P$ be a lattice polytope in $\rr^d$ and $p_1, \ldots, p_k$ an
ordered list 
of the lattice points in $P$.  The \emph{pulling triangulation}
$\Delta_{pull}(P)$ 
induced by this ordering is constructed recursively as follows:  If
$p_1, \ldots, p_k$ are affinely independent $\Delta_{pull}(P) = \{
\{p_1, \ldots, p_k\}\}$.  Else:
$$\Delta_{pull}(P) = \bigcup_F \{\{p_1\} \cup \sigma | \sigma \in
\Delta_{pull}(F) \} $$
where the union is over all facets $F$ of $P$ not containing $p_1$,
and the ordering of the lattice points in $F$ is the ordering induced
by the ordering of the lattice points in $P$.
\end{definition}

\begin{definition}
A triangulation $\Delta$ of a lattice polytope is called
\emph{unimodular} if every simplex in the triangulation attains the
minimal volume among all simplices formed by taking convex hulls of
lattice points in the polytope.  
\end{definition}

\begin{definition}
A lattice polytope $P$ is \emph{compressed} if every pulling
triangulation of $P$ using the lattice points in $P$ is unimodular.
If we are given a specific presentation of $P = P_A :=  \conv(A_1,
\ldots, A_n)$  as the convex hull of a finite set of integral points,
we say that $P_A$ is compressed if it 
is compressed with respect to the smallest lattice containing $A_1,
\ldots, A_n$.
\end{definition}

Compressed polytopes were introduced by 
Stanley in \cite{St} where unimodular was meant with respect to the
lattice $\zz^d$.  Our notion of unimodular is with respect to the
smallest lattice containing the integral points in $P$.  We say that
two lattice polytopes $P$ and $Q$ are \emph{lattice isomorphic} if there is
an affine isomorphism which is a bijection on their lattice points.
Our main
result is the following: 

\begin{theorem}\label{thm:main}
Let $\mathcal{L}$ be a lattice and suppose that $P$ is a lattice
polytope that has the irredundant 
linear description $P = \{ \xx \in \rr^d | a_i^T \xx \geq b_i, i = 1,
\ldots, n \}$.  Then the following conditions are equivalent:
\begin{enumerate}
\item $P$ is compressed.
\item For each $i$ there is at most one nonzero real number $m_i$ such
  that the set \\ $\{ \xx \in \mathcal{L} | a_i^T \xx = b_i + m_i\} \cap P$ is
  nonempty.
\item $P$ is lattice
  isomorphic to  an  integral polytope of the form $C_n \cap L$ where
  $C_n$ is the $n$-dimensional unit  hypercube and $L$ is an affine subspace.
\end{enumerate}
\end{theorem}

This result strengthens a result of Ohsugi and Hibi \cite{OH} who
essentially proved $(3) \implies (1)$.  Condition (2) suggests that
the term ``compressed'' is apt because compressed
polytopes are squeezed between two hyperplanes in every facet defining
direction. 

\begin{proof}
$(1) \implies (2)$  Supppose that $P$ is compressed and that for some 
$i$
there were two values $m > m'$ with  $\{ \xx \in \mathcal{L} | a_i^T \xx =
b_i + m\} \cap P$ and  $\{ \xx \in \mathcal{L} | a_i^T \xx = b_i + m'\} \cap P$
nonempty.  Let $p_m \in \{ \xx \in \mathcal{L} | a_i^T \xx =
b_i + m\} \cap P$ and $p_{m'} \in \{ \xx \in \mathcal{L} | a_i^T \xx = b_i +
m'\} \cap P$ and compare the pulling triangulations with $p_m$ first
and with $p_{m'}$ first and the same ordering of the lattice points in
the facet $F = \{ \xx \in \rr^d | a_i^T \xx = b_i\} \cap P$.  Then
given a simplex $\sigma$ in the pulling triangulation of $F$, the
ratio of volumes  ${\rm Vol}( p_m \cup \sigma)/ {\rm Vol}( p_n \cup
\sigma) = m/{m'} > 1$.  Hence the pulling triangulation of $P$ with $p_m$
first could not be unimodular contradicting the fact that $P$ was
compressed. 

$(2) \implies (3)$  Now suppose that $P$ satisfies condition (2)
above.  Since $P$ is a lattice polytope, condition (2) forces every
lattice point in $P$ to be a vertex since, given a facet defining
inequality $a^T \xx \geq b$, the largest value $m$ such that $P \cap \{
\xx \in \rr^d : a^T \xx = b + m \} $ is nonempty must have  $P \cap \{
\xx \in \mathcal{L} : a^T \xx = b + m \} $ nonempty as well as this set must
contain a vertex of $P$.  If there was a lattice point $p$ in $P$
which was not a vertex, it is in the relative interior of some face
$F$ of
$P$ of dimension greater than or equal to $1$.  This point $p$ could
not be in the set $P \cap \{\xx \in \rr^d : a^T \xx = b + m \}$ (where
$m$ is the unique largest value where this set is nonempty) for any
facet of $P$,  $a^T\xx = b$ which defines a nontrivial facet of $F$ and in
particular since $a^Tp > b$, there must be some value $m' < m$ such that
$P \cap \{\xx \in \mathcal{L} : a^T\xx = b + m' \}$ is nonempty. 

Now we must show that $P$ is affinely isomorphic to an integral
polytope that is the intersection of the unit hypercube with an affine
subspace.  Without  
loss of generality, we may suppose that $P$ does not lie in an affine
subspace: if 
it did we would make a unimodular change of coordinates to project to
a lower dimensional space.  This implies that in condition (2) there
is exactly 1 nonzero $m_i$ for each $i$. 
Consider the linear transformation $\pi:\rr^d \rightarrow \rr^n$
$$\xx \mapsto ( (a_1^T \xx - b_1)/m_1, \ldots, (a_n^T \xx - b_n)/m_n ).$$
The image $\pi(P)$ is a $0/1$ polytope since every vertex of $P$ is
mapped to a $0/1$ vector.  A point $p$ lies in $\pi(P)$ if and only if
$p \in C_n$ and $p$ is in the affine span of the image of the vertices
because the affine transformation $\pi$ sends the facet defining
inequality $a_i^T \xx \geq b_i$ to the inequality $y_i \geq 0$.  These
facts together imply that $P$ satisfies property $(3)$.

$(3) \implies (1)$  If the lattice polytope $P$ satisfies $(3)$ and
$\pi$ is the affine transformation, then $P$ is compressed if and only
if $Q = \pi(P)$ is compressed since this transformation maps the
lattice points in $P$ to the integer points in $Q$ and $P$ and $Q$ are
otherwise isomorphic.  Thus it remains to show that integral polytopes
$Q$ of the form $Q = C_n \cap \{\xx :A\xx = b \}$ are compressed.
This result is proven in \cite[Lemma 2.2]{OH}.  However, we will
provide a short self-contained proof of this fact.

Let $Q$ be an integral polytope of the form  $C_n \cap \{\xx :
A\xx = b \}$.  We will show $Q$ is compressed by induction on the
dimension.  If
$Q$ has dimension $0$ there is nothing to show.  Otherwise suppose $Q$
has dimension $d$ and consider any ordering of the vertices of $Q$.
Let $p$ be the first vertex and construct the pulling triangulation.
This is obtained by constructing the pulling triangulation of each
facet of $Q$ not containing $p$ and coning each of these
triangulations over $p$.  The normalized volume of each simplex is the
orthogonal distance from $p$ to the facet times the volume of
corresponding simplex in that facet.  However, each facet has
dimension $d-1$ and is 
of the form $C_n \cap \{ \xx : A\xx = b, x_i = 0 \}$ for some $i$ and
hence is compressed by induction.  Thus each simplex in the pulling
triangulation 
of each facet has normalized volume one.  Further, the orthogonal
distance to the 
corresponding facet is 1 since $p_i = 1$ when the facet is defined by
the equation $x_i = 0$.  So every simplex in the pulling triangulation is
unimodular.  Thus $Q$ is compressed.
\end{proof}

Many lattice polytopes which arise in applications (in particular, the
statistical applications from Section 5) possess symmetry groups that
are transitive on their lattice points.  From the preceding theorem we can
deduce that for such polytopes either every pulling triangulation is
unimodular or none are.

\begin{corollary}\label{cor:sym}
Suppose that $P$ is a lattice polytope and the group of affine
symmetries $\Gamma$ of $P$ is transitive on the lattice points of $P$.
Then either $P$ is compressed or no pulling triangulation of $P$ is
unimodular. 
\end{corollary}

\begin{proof}
We must show that if $P$ is not compressed then every pulling
triangulation is not unimodular.  To this end we can suppose that $P$
fails to satisfy condition (2) in the preceding theorem.  Then there
exists a facet $F = \{\xx \in \rr^d : a^T \xx = b \}$ of $P$ and two
nonzero reals $m > m'$ such that $\{ \xx \in \mathcal{L} : a^T \xx = b + m'\}
\cap P$ and  $\{ \xx \in \mathcal{L} : a^T \xx = b + m\}
\cap P$ are nonempty.  Consider any ordering of the vertices of $P$
and the resulting pulling triangulation.  After applying a suitable
element $g \in \Gamma$ to this ordering, we can assume that the first
point $p_m$ in the 
pulling triangulation is in the set $\{ \xx \in \mathcal{L} : a^T \xx = b + m\}
\cap P$.  Consider any other pulling triangulation which
has the same order of the points in $F$ and a point $p_{m'}$ in $\{ \xx \in
\mathcal{L} : a^T \xx = b + m'\} \cap P$ as the first vertex.   Among the
simplices in the first pulling triangulation of $P$ are those of the
form $p_m \cup \sigma$ and in the second pulling triangulation $p_{m'}
\cup \sigma$ where $\sigma$ is in the induced pulling triangulation of
$F$.  We see that the ratio of volumes of these simplices  $Vol( p_m \cup
\sigma) / Vol(p_{m'} \cup \sigma) = m/{m'}$ and hence the first pulling
triangulation could not be unimodular.  However, this pulling
triangulation was arbitrary, so no  pulling triangulation of $P$ is
unimodular. 
\end{proof}

\section{Compressed cut polytopes}

As an application of our characterization of compressed polytopes, we
describe those graphs  $G$ whose cut polytopes are compressed.
We assume throughout that $G = (V_n,E)$ is an undirected graph with
vertices $V_n = [n] := \{1,2,\ldots, n \}$ and edges $E$ without loops
or multiple edges.  Our definitions and notation comes from \cite{DL}
and we assume some familiarity with the basic facts about these polytopes.

\begin{definition}
Let $S \subseteq V_n$.  The cut semimetric on $G$ induced by $S$
is the 0/1 vector $\delta_{G}(S)$ in $\rr^E$ defined by
$$\delta_G(S)_{ij} = 1 \mbox{ if } |S \cap \{i,j\}| = 1, \mbox{ and }
\delta_G(S)_{ij} = 0 \mbox{ otherwise},$$
where $ij \in E$.  The cut polytope of $G$ is the 0/1 polytope
$$\cut(G) = \conv( \delta_G(S) | S \subseteq V_n ).$$
\end{definition}

We will apply criterion (2) from the main theorem to deduce the
following:

\begin{theorem}\label{thm:cut}
The cut polytope $\cut(G)$ of a graph $G$ is compressed if and only if $G$ has
no $K_5$ minors and every induced cycle in $G$ has length less than or equal
to $4$. 
\end{theorem}

A cycle in a graph is induced if there is no chord in the graph
cutting across it.  Equivalently, a cycle is induced if it is an
induced subgraph.  The proof of the theorem requires a few
intermediate results.

\begin{lemma}
If $\cut(G)$ is compressed and $H$ is obtained from $G$ by contracting
an edge then $\cut(H)$ is compressed. 
\end{lemma}

\begin{proof}
Let $ij$ be the contracted edge.  The polytope $\cut(H)$ is isomorphic
to $\{\xx | x_{ij} = 0 \} \cap 
\cut(G)$ and hence is isomorphic to a face of $\cut(G)$.  But every
face of a compressed polytope is compressed. 
\end{proof}

\begin{lemma}\label{lem:subgraph}
If $\cut(G)$ is compressed and $H$ is an induced subgraph of $G$ then
$\cut(H)$ is compressed.
\end{lemma}

\begin{proof}
Let $E' \subset E$ be the union of all edges in $G$ not incident to
$H$ together with exactly one edge which is incident to $H$ but not
contained in $H$ (provided such an edge exists).  Then $\cut(H)$ is
isomorphic to $\{\xx | x_e = 0, e \in E' \} \cap \cut(G)$, and hence
is isomorphic to a face of $\cut(G)$.  But every face of a compressed
polytope is compressed.
\end{proof}

\begin{lemma}
The polytope $\cut(K_5)$ is not compressed.
\end{lemma}

\begin{proof}
One facet defining inequality for $\cut(K_5)$ comes by via the
following hypermetric construction \cite{DL}.  Let $b = (1,1,1,-1,-1)$
and consider the inequality
$$\sum_{1 \leq i < j \leq 5} b_ib_j x_{ij} \leq 0.$$
This inequality defines a facet of $\cut(K_5)$ called a pentagonal
facet.  To show that $\cut(K_5)$ is not compressed it suffices to
exhibit two sets $S, T \subset V_5$ such that

$$\sum_{1 \leq i < j \leq 5} b_ib_j \delta_{K_5}(S)_{ij} < \sum_{1
  \leq i < j \leq 5} b_ib_j \delta_{K_5}(T)_{ij} < 0,$$ 

\noindent since the cut semimetrics are integral points in the cut
polytope.  Taking $S = \{1,2,3\}$ and $T = \{1,2\}$ yields

$$ -6 = \sum_{1 \leq i < j \leq 5} b_ib_j \delta_{K_5}(S)_{ij} < \sum_{1
  \leq i < j \leq 5} b_ib_j \delta_{K_5}(T)_{ij}= -2 < 0.$$ 
 
\end{proof}

The preceding three Lemmas imply that if we want to identify graphs
whose cut polytopes are compressed, we may restrict attention to those
graphs without $K_5$ minors.  In general, it remains a hard open
problem to give a facet description of the cut polytopes, however, in
the special case of graphs without $K_5$ minors, a complete
irredundant linear description is known.

\begin{theorem}\label{thm:k5}
Let $G$ be a graph without $K_5$ minors.  Then $\cut(G)$ is the
solution set of the following linear inequalities:
$$0 \leq x_{e} \leq 1, \, e \in E$$
$$\sum_{e \in F} x_e - \sum_{e \in C \setminus F} x_e \leq |F| - 1$$
where $C$ ranges over the induced cycles of $G$ and $F$ ranges over
the odd subsets of $C$.  Each of the linear inequalities of the
second type is facet defining and the inequalities $0 \leq x_e \leq 1$
may or may not be facet defining.
\end{theorem}

Theorem \ref{thm:k5} is a consequence of the decomposition theory for
binary matroids.  It is proven in \cite{BM} and depends on results in
\cite{Sey}.
Thus to prove the main theorem in this section we just need to
determine under what conditions these facet defining inequalities
satisfy condition $(2)$ from Theorem \ref{thm:main}.  For the
inequalities of type $0 \leq x_e \leq 1$, these always satisfy
condition (2) regardless of whether or not they are facet defining.
Since the 
structure of the remaining facet defining inequalities only depends on
the induced cycles  in the graph it suffices to prove the following:

\begin{lemma}
Let be $C$ an induced cycle of $G$, and
$F$ an odd subset of $C$.  Then the set
$$\{\xx \in \zz^d | \sum_{e \in F}  x_e - \sum_{e \in C \setminus F}
  x_e = |F| - 1 - m \} \cap \cut(G)$$
is nonempty for exactly $\lfloor\frac{|C|}{2}\rfloor - 1$ nonzero
values of $m$.  
\end{lemma}

\begin{proof}
Since the value of the linear functional $\sum_{e \in F}  x_e - \sum_{e
\in C \setminus F} x_e$ only depends on the edges in $C$, we can
assume that $G = C$.  Furthermore, the operation of switching (see
\cite{DL}) shows that each such facet is equivalent (i.e. up to change
of coordinates) to the facet given
by $x_{12} - x_{23} - \ldots - x_{1n} \leq 0$.  So it suffices to
prove the Lemma in this setting.   

Since cut semimetrics $\delta_G(S)$ are the only integral points in
$\cut(G)$ it suffices to determine what values $\delta_G(S)_{12} -
\delta_G(S)_{23} - \ldots - \delta_G(S)_{1n}$ can take.  Modulo 2, 
$$\delta_G(S)_{12} -
\delta_G(S)_{23} - \cdots - \delta_G(S)_{1n} = \delta_G(S)_{12} +
\delta_G(S)_{23} + \cdots +\delta_G(S)_{1n} \equiv 0 \mod 2$$

\noindent so $\delta_G(S)_{12} -
\delta_G(S)_{23} - \ldots - \delta_G(S)_{1n}$ must be even.  Since
$\delta_G(S)_{ij}$ is either a zero or a one, there are at most
$\lfloor\frac{|C|}{2}\rfloor - 1$ nonzero values that this expression can
take.  However, for each $j$ with $0 < j \leq \lfloor\frac{|C|}{2}\rfloor$ the
set $S_j = \{ 2i : i \in [j] \}$ has 
$$\delta_G(S_j)_{12} -
\delta_G(S_j)_{23} - \cdots - \delta_G(S_j)_{1n} = 2-2j$$
which completes the proof.
\end{proof}

\section{Compressed polytopes in linear optimization}

Compressed polytopes are closedly tied to linear integer
optimization problems.  In particular, we consider the following setup.
Let $A$ be an integral matrix with columns $A_1, A_2, \ldots, A_n$.
We assume throughout that $A$ is homogeneous in the sense that there
is a nonzero weight vector $w$ such that $w^T A_i = 1$ for all $i$.
For each $i$ consider the integer programming problem

$$\mbox{Maximize } x_i \mbox{ subject to}$$
$$A \xx = b, \xx \geq 0, \xx \mbox{ integral.}$$

\noindent For a given $i$, $A$ and $b$ we denote the optimal value of the
integer program by $IP^+_i(A,b)$.  We call a vector $b$ IP-feasible if
$b = A \xx$ for some nonnegative integral $\xx$.  The corresponding linear
programming relaxation drops the integrality consideration:

$$\mbox{Maximize } x_i \mbox{ subject to}$$
$$A \xx = b, \xx \geq 0.$$

\noindent We denote the optimal value of the linear programming relaxation by
$LP^+_i(A,b)$.  Since linear programs are considerably easier to solve
than integer programs, a fundamental question in optimization is to decide
what conditions guarantee that $LP^+_i(A,b) = IP^+_i(A,b)$.
Let $P_A$ be the polytope obtained by taking the convex 
hull of the columns of $A$.  Pulling triangulations of $P_A$ provide a useful
sufficient condition to guarantee $LP^+_i(A,b) = IP^+_i(A,b)$.

\begin{proposition}\label{thm:pullip}
For fixed $A$ and $i$, $LP^+_i(A,b) = IP^+_i(A,b)$ for all IP-feasible $b$ if
there exists some  
ordering of the columns of $A$ with $A_i$ first such that the
pulling triangulation of $P_A$ using only $A_1, \ldots, A_n$ is
unimodular. 
\end{proposition}

\begin{proof}
We provide a sketch of the proof which depends on some well known
results in computational algebra.  Details can be found in
\cite[Chapter 8]{St}.
The linear programming relaxation solves the standard form integer
program for all right hand sides $b$ if an associated
intial ideal of the toric ideal $I_A$ is squarefree.  The initial
ideal is squarefree if and only if the corresponding regular
triangulation of $P_A$ is unimodular.  In the case where the
associated cost vector is the maximization of the $x_{i}$ coordinate, the
corresponding triangulation is a pulling triangulation of $P_A$ with
$A_{i}$ first.   
\end{proof}

The condition in Proposition \ref{thm:pullip} is not, however,
necessary:  if $LP^+_i(A,b) = IP^+_i(A,b)$ for all $b$ there need not
exist a unimodular pulling triangulation of $P_A$ with $A_i$ first as
the following example illustrates.

\begin{example}
Consider the matrix $A$ given by

$$A = \begin{pmatrix}
1 & 1 & 1 & 1 & 1 \\
0 & 0 & 1 & 2 & 3 \\
1 & 0 & 0 & 0 & 0 \end{pmatrix}.$$

This matrix has the property that $LP^+_1(A,b) = IP^+_1(A,b)$ for all
IP-freasible $b$.  Indeed, given an $IP$ feasible $b$, every
nonnegative vector $\xx$ with $A \xx = b$ has $x_1 = b_3$.  On the
other hand, $P_A$ has no unimodular pulling triangulations.
\end{example}  

This subtlety drops away if we require that $LP^+_i(A,b) =
IP^+_i(A,b)$ for all IP-feasible $b$ \emph{and} for all $i$.

\begin{theorem}\label{thm:opt}
Let $A$ be a homogeneous matrix.  Then $LP^+_i(A,b) = IP^+_i(A,b)$ for all
$i$ and all IP-feasible $b$ if and only if $P_A$ is compressed.
\end{theorem}

Recall that in this context where $P_A = \conv (A_1, \ldots, A_n)$ we
mean that $P_A$ is compressed with respect to the largest lattice
containing $A_1, \ldots, A_n$.

\begin{proof}
If $P_A$ is compressed then any pulling triangulation with $A_i$
first is unimodular which implies by Proposition \ref{thm:pullip} that the
LP optimums equal the IP optimums.  
Conversely, if $P_A$ is not
compressed, there is a facet 
defining inequality which violates condition (2) in the main theorem.
We will use this violation to construct an IP feasible $b$ such that
the LP optimum for the maximization problem cannot equal the IP
optimum.  

Denote the violating facet by $F = \{ \xx \in \rr^d: a^T \xx = b \}$.  Since
$P_A$ is a polytope, there is a largest real number $m$ such that $
\{\xx \in \mathcal{L} : a^T \xx = b + m \}$  is nonempty.  We may suppose that $A_1 \in
\{\xx \in \mathcal{L} :a^T \xx = b+ m \}$.  We will partition $A_2, \ldots, A_n$ in
the following manner: $a^T A_i = b+m$ for $i = 2, \ldots, k$, $b <
a^T A_i < b+m $ for $i = k+1, \ldots, l$, and $a^T A_i = b$ for $i =
l+1, \ldots, n$.  Let 
$$K  = \ker_\zz(A) \cap  \{{\bf y} \in \zz^d | y_1 < 0, 
y_2 \leq 0, \ldots, y_k \leq 0, y_{k+1} \geq 0,y_{k+2} \leq 0 \ldots, y_l
\leq 0\}.$$

Note that $K$ is nonempty.  This follows since there exist
affine dependencies among $A_1$, the elements of $F \cap \mathcal{L}$,
and $A_{k+1}$ (there are at least $d+2$
points in a $d$ dimensional lattice).  Furthermore, any such affine
dependency must have $y_1$ and $y_{k+1}$ with opposite signs since neither
$A_1$ nor $A_{k+1}$ are contained in $F$.  Among the vectors in $K$, let
$\vv \in K$ be any such vector with $v_{k+1}$ with the
minimal value among all $\vv$ in $K$.  This minimal value is
strictly greater than $1$.  We define the right-hand side vector $b$ which will
violate $LP^+_1(A,b) > IP^+_1(A,b)$ by 
$$b = \sum_{i | v_i > 0 } v_i A_i - A_{k+1}.$$ 
Clearly $b$ is IP-feasible since we have expressed it as a nonnegative
combination of the columns of $A$.  

First of all, we claim that $IP^+_1(A,b) = 0$.  If not, there is an improving
integer vector $\uu \in K$ with such that the vector $\vv^+ - e_{k+1} -\uu$
is nonegative and has first coordinate greater than zero.  The
existence of such a $\uu$ violates our minimality assumption on
$v_{k+1}$ (since $u_{k+1} \leq v_{k+1} - 1$).  On the 
other hand, the rational vector $\uu = \frac{v_{k+1}-1}{v_{k+1}} \vv$ is
an improving vector such that $\tilde{\vv} = \vv^+ - e_{k+1} - \uu$ is a nonnegative
rational vector with $A\tilde{\vv} = b$ and $\tilde{v}_1 > 0$ so that
$LP^+_1(A,b) > 0$.      
\end{proof}

\section{Applications in statistical disclosure limitation}

One motivation for studying compressed polytopes comes from their
relationship to certain optimization problems which arise in
statistical disclosure limitation.  The general problem in this area
is to determine what information about individual survey respondents
can be inferred from the release of partial data.  This type of
problem arises when government agencies like a census bureau gather
information about citizens and wish to release partial data to the
public for the purposes of data analysis but are required by law to maintain
the privacy of citizens.

The case we consider here concerns the release of margins of a
multiway contingency table.  In this case, an individual cell entry is
considered secure if among all nonnegative integral tables with given
released marginal totals the upper and lower bounds on the cell
entry are far enough apart \cite{BG,Ch}.  This naturally leads to
standard form integer programs of the following type:

$$\mbox{Maximize/ Minimize } x_{\bf 0} \mbox{ subject to}$$
$$ A_\Delta \xx = b, \xx \geq 0, \mbox{ and } \xx \mbox{ integral},$$

\noindent where $A_\Delta$ is a certain $0/1$ matrix which computes
the released 
margins $b$ of the multiway table $\xx$.  A heuristic for
approximating the solution to this integer program is to solve the
linear programming relaxations:

$$\mbox{Maximize/ Minimize } x_{\bf 0} \mbox{ subject to}$$
$$ A_\Delta\xx = b \mbox{ and }  \xx \geq 0.$$

A fundamental problem in this area is to determine under what
conditions the linear programming relaxation is equal to the true
integer value. We will focus here on the maximization problem.  
To state our results, we first need to establish notation
for the contingency table problems of interest.  Here $\xx$ denotes a
$d_1 \times d_2 \times \cdots \times d_n$ multiway contingency table.
The particular collection of margins of this table which are released
are encoded by a simplicial complex $\Delta$ on the $n$-element set
$[n]$.  

Each facet $S \in \Delta$ corresponds to a released margin.
Computing a collection of marginals of a multiway table is a linear
transformation.  The matrix, represented in the standard basis, which
encodes this linear transformation is denoted by $A_\Delta$.   Note
that the size of the matrix $A_\Delta$ and problems related to linear
programming relaxations depend on $\Delta$ \emph{and} the integer vector $d =
(d_1, d_2, \ldots, d_n)$ though we suppress the dependence on $d$ when we
use the notation $A_\Delta$.  We use the notation $P_\Delta$ to denote
the convex hull of the columns of the matrix $A_\Delta$.  From the
previous section, we deduce the following basic fact:

\begin{corollary}
The linear programming relaxation solves the integer programs
$IP^+_{\bf 0}(A_\Delta, b) = LP^+_{\bf 0}(A_\Delta,b)$ for all
marginals $b$ if and only if the marginal polytope $P_\Delta$ is compressed.  
\end{corollary}

\begin{proof}
Becuase of the transitive symmetry group on the vertices of
$P_\Delta$, $IP^+_{\bf 0}(A_\Delta,b) = 
LP^+_{\bf 0}(A_\Delta,b)$ for the ${\bf 0}$ cell entry implies if and
only if
this holds for all cell entries.  Then by Theorem \ref{thm:opt} this
holds if and only if $P_\Delta$ is compressed.
\end{proof}

Thus we are led to study the following general problem:  

\begin{problem}
Characterize the pairs $(\Delta, d)$ of simplicial complexes $\Delta$
and integer vectors $d = (d_1, \ldots, d_n)$ such that $P_\Delta$ is
compressed. 
\end{problem}

It seems a
challenging problem to classify such marginals in general, since it
would require the knowledge of many families of facet defining
inequalities of the marginal polytopes.  There is very little known
about these facet defining inequalities in general.  In the
remainder of this section, we provide some constructions for producing
compressed marginal polytopes.  As a corollary, we deduce that the
marginal polytopes of decomposable models are compressed.
We also provide a complete
characterization of compressed marginal polytopes in two restricted cases.

There are a few standard operations
on simplicial complexes that send compressed
marginal polytopes to compressed  marginal polytopes.

\begin{proposition}\label{prop:down}
Suppose that the pair $(\Delta, d)$ has $P_\Delta$ compressed.
\begin{enumerate}
\item If $\Delta' \subset
\Delta$ is an induced subcomplex and $d'$ the correspond integer
vector then the pair $(\Delta',d')$ has $P_{\Delta'}$ compressed.
\item If $d' \leq d$ coordinatate-wise then the pair $(\Delta',d')$
  with $\Delta' = \Delta$ has $P_{\Delta'}$ compressed.
\end{enumerate}
\end{proposition}

\begin{proof}
In both cases $P_{\Delta'}$ is isomorphic to a face of $P_\Delta$.
However, the faces of compressed polytopes are compressed.
\end{proof}

\begin{proposition}
Suppose that the pair $(\Delta, d)$ has the marginal polytope
$P_\Delta$ compressed.  And let $\Delta'$ be the new simplicial
complex on $[n+1]$ obtained from $\Delta$ by
$\Delta' = \{ \{n+1\} \cup F | F \in \Delta \}$
and $d' = (d_1, \ldots, d_n, d_{n+1})$ where $d_{n+1}$ is any positive
integer.  Then the pair $(\Delta',d')$ has a compressed marginal
polytope $P_{\Delta'}$.  
\end{proposition}

\begin{proof}
The marginal polytope $P_{\Delta'}$ is isomorphic to the direct join
of $d_{n+1}$ copies of $P_{\Delta}$.  But the direct join of compressed
polytopes is compressed since any triangulation of the direct join is
obtained by taking the direct join of the induced triangulations of
the pieces.  The direct join of two unimodular triangulations is unimodular.
\end{proof}

\begin{definition}
A simplicial complex $\Delta$ is called \emph{reducible} with decomposition
$(\Delta_1, S, \Delta_2)$ if
\begin{enumerate}
\item $\Delta_1$ and $\Delta_2$ are induced subcomplexes of $\Delta$,
\item $S \subset [n]$,
\item $\Delta_1 \cup \Delta_2 = \Delta$, and
\item $\Delta_1 \cap \Delta_2 = 2^S$.
\end{enumerate}
A simplicial complex is called \emph{decomposable} if $\Delta$ is reducible
and each of $\Delta_1$ and $\Delta_2$ is either decomposable or a simplex.
\end{definition}

Given a reducible simplicial complex $\Delta$ with decomposition
$(\Delta_1, S, \Delta_2)$ together with the integer vector $d$ denote
by $d^1$ and $d^2$ the induced vectors with indices corresponding to
the nodes of $\Delta_1$ and $\Delta_2$ respectively.

\begin{proposition}\label{prop:red}
If $\Delta$ is reducible and the pairs $(\Delta_1, d^1)$ and
$(\Delta_2,d^2)$ have compressed marginal polytopes then the marginal
polytope $P_\Delta$ is compressed.
\end{proposition}

\begin{proof}
For reducible models $\Delta$, the marginal polytopes are given by
$$P_\Delta = P_{\Delta_1} \times P_{\Delta_2} \cap \{(\xx,\yy) |
\pi_1(\xx) = \pi_2(\yy) \}$$
where $\pi_1$ and $\pi_2$ are the $S$-marginal maps of $\xx$ and $\yy$
repsectively.  In particular, the set of facet defining inequalities
of $P_\Delta$ is just the union of the facet defining of
$P_{\Delta_1}$ and $P_{\Delta_2}$.  Since $P_{\Delta_1}$ and
$P_{\Delta_2}$ are compressed, these facet defining inequalites
satisfy condition (2) of the main theorem.  But this implies that they
also satisfay condition (2) of the main theorem with respect to
$P_{\Delta}$ as well.  This implies that $P_\Delta$ is compressed.
\end{proof}

\begin{corollary}
If $\Delta$ is decomposable then $P_\Delta$ is compressed.
\end{corollary}

\begin{proof}
If $\Delta = 2^{[n]}$ then $P_\Delta$ is a simplex.  Thus, if $\Delta$
is decomposable $P_\Delta$ is compressed by applying Proposition
\ref{prop:red} and induction on the number of facets of $\Delta$.
\end{proof}

The preceding propositions provide methods for producing compressed
marginal polytopes from smaller compressed marginal polytopes, however
these results are far from giving a complete characterization of all
pairs $(\Delta,d)$ such the marginal polytopes are compressed.
In the remainder of this section, we
provide characterizations of compressed marginal polytopes in two
settings where we place ``extremal'' conditions on $\Delta$, or $d$ or both.

\begin{proposition}
Let $\Delta$ be the boundary of an $n-1$ simplex.  Then $P_\Delta$ is
compressed if and only if for at most two $i$, $d_i > 2$ or $n = 3$
and up to symmetry $d = (3,3,d_3)$.
\end{proposition}

\begin{proof}
In the case where for at most two $i$, $d_i >2$, it is known that
$P_\Delta$ is a unimodular polytope (e.g. \cite[Chapter 14]{Sturm})
and hence is compressed.
The case where $n = 3$ and $d = (3,3,d_3)$, the complete facet
description of this polytope is known (e.g. \cite{EFRS}) and one
verifies 
that the facet defining inequalities in this case satisfy condition
(2) in the main theorem.  Direct computation shows that condition (2)
of the main theorem fails in the case $n
= 3$, $d = (3,4,4)$ and $n = 4$, $d = (2,3,3,3)$.  These results
together with Proposition \ref{prop:down} imply that $P_\Delta$ is
compressed in no other cases.
\end{proof}

The cut polytopes from the previous section are intimately tied to the
marginal polytopes we are interested in, in the special case where $d
= (2,2, \ldots, 2)$ and all facets of $\Delta$ are 0 or
1-dimensional.  In this case $\Delta$ is a graph and we have the
following well known result (see \cite{DL}):

\begin{lemma}
Given a graph $\Delta$ and $d = (2,2,\ldots, 2)$ there is an affine
isomorphism  of the marginal polytope $P_\Delta$ to the cut polytope
$\cut(\widetilde{\Delta})$ where $\widetilde{\Delta}$ is the graph
obtained from $\Delta$ by adding a new vertex $v$ and all edges from
$v$ to the nodes of $\Delta$.
\end{lemma}

The affine isomorphism in the preceding Lemma is known as the
covariance mapping.  Then we can deduce:

\begin{theorem}
Let $\Delta$ be a graph and $d = (2,2, \ldots, 2)$.  Then $P_\Delta$ is compressed  if and only if
$\Delta$ is free of $K_4$ minors and every induced
cycle in $\Delta$ has length less than or equal to $4$.  
\end{theorem}

\begin{proof}
The graph $\widetilde{\Delta}$ is free of $K_5$ minors and has all
induced cycles of length less than or equal to four if and only if
$\Delta$ has no $K_4$ minors and all induced cycles of length less
than or equal to four.  Thus this is a direct consequence of Theorem
\ref{thm:cut} which characterized the compressed cut polytopes.
\end{proof}

In these cases we can in fact say more:  even though the size of the
integer program seems exponential in $n$ the number of nodes in the
simplicial complex, in the case where $P_\Delta$ is compressed we can
solve the 
corresponding linear program (and hence the integer program) in
polynomial time. 

\begin{corollary}
Suppose that $d = (2,2, \ldots, 2)$ and $\Delta$ is a graph that is
free of $K_4$ minors and has every induced cycle of length less than
or equal to four. Then the IP-maximum value $IP^+_i(A_\Delta,b)$ can
be computed in polynomial time in $n$ and the bit complexity of $b$. 
\end{corollary}

\begin{proof}
Since $IP^+_i(A,b) = LP^+_i(A,b)$ for these graphs, it suffices to show
that the linear program can be solved in polynomial time.  However,
the problem of maximizing a coordinate is polynomial time equivalent
to determining if a point lies in $P_\Delta$.  For graphs without $K_4$
minors, the containment problem can be decided in polynomial time as
illustrated in \cite{DL}.
\end{proof}

In general, we would like to understand how far the linear programming relaxations can be from the true integer programming values for these optimization problems in statistical disclosure limitation.  This leads to the study of the integer programming gap \cite{HS}.  A natural question to ask is: How does the failure of condition (2) in Theorem \ref{thm:main} relate to the integer programming gap?  A natural family of marginal polytopes where this problem could be explored is the family of cycles.


\bigskip
\address{ 
Department of Mathematics \\
University of California \\
Berkeley, CA 94720-3840 \\
USA
}
{seths@math.berkeley.edu}


\begin{thebibliography}{99}

\bibitem{BM} 
F.~Barahona and A.~R.~Mahjoub.  On the cut
polytope. \emph{Mathematical Programming}, {\bf 36}: 157--173, 1986.

\bibitem{BG}
L. Buzzigoli and A. Giusti.  An algorithm to calculate the lower
and upper bounds of the elements of an array given its marginals,
in \emph{Statistical Data Protection Proceedings}, Eurostat,
Luxembourg (1999), 131--147.

\bibitem{Ch}
S.~D.~Chowdhury, G.~T.~Duncan, R.~Krishnan, S.~F.~Roehrig  and
S.~Mukherjee. Disclosure Detection in Multivariate 
Categorical Databases: Auditing Confidentiality Protection Through
Two New Matrix Operators. {\sl Management Science} (1999) {\bf 45} No.
12, 1710--23.

\bibitem{DL}
M.~M.~Deza and M.~Laurent. \emph{Geometry of Cuts and
  Metrics}. Algorithms and Combinatorics {\bf 15}  Springer-Verlag,
Berlin, 1997.

\bibitem{DS} P.~Diaconis and B.~Sturmfels.
Algebraic algorithms for sampling from conditional distributions.
\emph{Annals of Statistics}, 26 (1998), 363--397.

\bibitem{EFRS} S.~E.~Fienberg, N.~Eriksson, A.~Rinaldo, and S.~Sullivant.
Polyhedral conditions for the nonexistence of the MLE for hierarchical
log-linear models.
To appear in Journal of Symbolic Computation, Special issue on
Computational Algebraic Statistics,  {\tt math.CO/0405044}, 2004. 

\bibitem{HS}
S. Ho\c{s}ten and B. Sturmfels.  Computing the integer
  programming gap.  To appear in \emph{Combinatorica}, 2003.

\bibitem{OH}
H.~Ohsugi and T.~Hibi.  Convex polytopes all of whose reverse
lexicographic initial ideals are squarefree.
\emph{Proc. Amer. Math. Soc.} {\bf 129}   (2001), 2541--2546.  

\bibitem{Sey}
P.~D.~Seymour.  Matroids and multicommodity flows.  \emph{European
  Journal of Combinatorics}. {\bf 2}:257--290, 1981.

\bibitem{St}
R.~Stanley.  Decompositions of rational convex polytopes.
\emph{Ann. Discrete Math.} {\bf  6} (1980), 333 -- 342. 

\bibitem{Sturm}
B. Sturmfels.  \emph{Gr\"obner Bases and Convex Polytopes},
American Mathematical Society. Providence, RI, 1995.

\bibitem{Z}
G. Ziegler.  \emph{Lectures on Polytopes}. Graduate Texts in
Mathematics.  Springer-Verlag, New York, 1995.


\end{thebibliography}
\end{document}